\documentclass[12pt, reqno, a4paper]{amsart}
\usepackage{csquotes}
\usepackage{graphicx,amsaddr}
\usepackage{amsfonts}
\usepackage{amsmath}
\usepackage{mathtools}
\usepackage{amsthm}
\usepackage{amsfonts}
\usepackage{amssymb}
\usepackage{calc}
\usepackage{rotating}
\usepackage[all,cmtip]{xy}
\usepackage{pifont}
\usepackage{enumitem}
\usepackage{lscape}
\usepackage{array}

\theoremstyle{plain}
\newtheorem{defn}{Definition}[section]
\newtheorem{thm}[defn]{Theorem}
\newtheorem{coro}[defn]{Corollary}

\newtheorem{exam}[defn]{Example}

\numberwithin{equation}{section}

\title[Spatial Numerical Range]{Spatial Numerical Range in Non-unital, Normed algebras and their Unitizations}

\author{H. V. Dedania}

\address{Dept. of Mathematics, Sardar Patel University, Vallabh Vidyanagar 388120, Gujarat, India}
\email{hvdedania@gmail.com}

\author{A. B. Patel*}

\address{Dept. of Mathematics, Sardar Patel University, Vallabh Vidyanagar 388120, Gujarat, India}
\email{avadhpatel663@gmail.com}
\newcolumntype{C}[1]{>{\centering\arraybackslash}m{#1}}

\begin{document}



\maketitle
\section*{Abstract}
Let $(A, \|\cdot\|)$ be any normed algebra (not necessarily complete nor unital). Let $a \in A$ and let $V_A(a)$ denote the spatial numerical range of $a$ in $(A, \|\cdot\|)$. Let $A_e = A + {\mathbb C} 1$ be the unitization of $A$. If $A$ is faithful, then we get two norms on $A_e$; namely, the operator norm $\|\cdot\|_{op}$ and the $\ell^1$-norm $\|\cdot\|_1$. Let $A^{op} = (A, \|\cdot\|_{op})$, $A_e^{op} = (A_e, \|\cdot\|_{op})$, and $A_e^1 = (A_e, \|\cdot\|_1)$. We can calculate the spatial numerical range of $a$ in all these three normed algebras. Because the spatial numerical range highly depend on the identity as well as on the completeness and the regularity of the norm, they are different. In this paper, we study the relations among them. Most of the results proved in \cite{BoDu:71, BoDu:73} will become corollaries of our results. We shall also show that the completeness and regularity of the norm is not required in \cite[Theorem 2.3]{GaHu:89}.

\section{Introduction }\label{intro}
Throughout $A$ is any normed algebra; the (algebra) norm on $A$ will be denoted by $\|\cdot\|$. The algebra $A$ is \emph{faithful} if $a=0$ whenever $a \in A$ and $aA = \{0\}$. Let $A_e = A + {\mathbb C}1$ be the unitization of $A$~\cite[page. 15]{BoDu:73}, where $1$ is the identity of $A_e$. If $A$ is faithful, then any norm $\|\cdot\|$ on $A$ induces following two norms on $A_e$; namely, the operator norm $\|\cdot\|_{op}$ and the $\ell^1$-norm $\|\cdot\|_1$, which are defined as follow.
\begin{eqnarray*}
  \|a+\lambda 1\|_{op} &=& \sup\{\|a x + \lambda x \|: x \in A, \|x\| \leq 1\} \\
  \|a+\lambda 1\|_1 &=& \|a\| + |\lambda|.
\end{eqnarray*}
In general, $\|\cdot\|_{op} \leq \|\cdot\|$ on $A$. The norm $\|\cdot\|$ is \emph{regular} if $\|\cdot\|_{op} = \|\cdot\|$ on $A$. We shall use the notions $A = (A, \|\cdot\|)$, $A^{op} = (A, \|\cdot\|_{op})$, $A_e^{op} = (A_e, \|\cdot\|_{op})$, and $A_e^1 = (A_e, \|\cdot\|_1)$.

Let $S(A)= \{x \in A: \|x\|=1\}$ be the unit sphere in $A$. Let $A^*$ be the Banach space dual of $A$. Let $D_{A}(x) = \{\varphi \in A^*: \|\varphi\|=1=\varphi(x)\}$ for each $x \in S(A)$. Further, let $V_{A}(a;x) = \{\varphi(ax): \varphi \in D_{A}(x)\}$ for $a \in A$ and $x \in S(A)$. Then $V_{A}(a) = \cup\{V_{A}(a;x): x \in S(A)\}$ is the \emph{spatial numerical range (SNR)} and $\nu_{A}(a) = \sup \{|\lambda|: \lambda \in V_{A}(a)\}$ is the \emph{spatial numerical radius} of $a$ in $(A, \|\cdot\|)$. It is still open problem wether $V_A(a)$ is always convex? This problem is discussed in \cite{DePa:23}.

We shall see that the spatial numerical range $V_{A}(a)$  highly depends on both the algebra $A$ and the norm $\|\cdot\|$. Therefore the sets $V_{A}(a)$, $V_{A^{op}}(a)$, $V_{A_e^{op}}(a)$, and $V_{A_e^1}(a)$ cannot be identical. In this paper, we study the relations among them and exhibit various examples. We have applied this concept in proving some results on the spectral extension property (SEP) in non-unital Banach algebras~\cite{DePa:22}.

\section{Main Results}\label{sec:2}

Recall that the spatial numerical range $V_A(a) = \cup\{V_A(a;x) : x \in S(A)\}$. So first we list out some basic properties of the sets $V_A(a;x)$.

\begin{thm}\label{r2}
Let A be any normed algebra, $a,b \in A, \, x,y \in S(A)$, and $\alpha \in \mathbb C$. Then
\begin{enumerate}
  \item $V_{A}(\alpha a;x)=\alpha V_{A}(a;x)$ and $V_{A}(a+b;x) \subset V_{A}(a;x)+ V_{A}(b;x)$;
  \item $V_A(a;x)$ is compact and convex;
  \item $V_{A}(a;x) = V_{A}(a;y)$ whenever $x$ and $y$ are linearly dependent;
  \item Let $B$ be a subalgebra of $A$. Then $$V_B(a;x) = V_A(a;x)  \; (a \in B, \, x \in S(B));$$
  \item Let $a_n \longrightarrow a$ in $A$. Then $\overline{V_A(a;x)} \subset \overline{\cup_{n=1}^{\infty}V_A(a_n;x)}$.\label{r2a}
\end{enumerate}
\end{thm}

\begin{proof}
(1) This is easy.

(2) Let $\{\lambda_n\} \subset V_A(a;x)$ such that $\lambda_n \longrightarrow \lambda$. Choose $x \in S(A)$ and $\varphi_n \in D_A(x)$ such that $\varphi_n(ax)=\lambda_n$. Since $D_A(x)$ is weak*-compact, we may assume that $\varphi_n \longrightarrow \varphi$ in $D_A(x)$ in weak*-topology. Then $\lambda_n=\varphi_n(ax) \longrightarrow \varphi(ax)$. Hence $\lambda = \varphi(ax) \in V_A(a;x)$. So $V_A(a;x)$ is closed. Clearly, it is a bounded set.
Let $\lambda_1, \lambda_2 \in V_{A}(a;x)$. Then there exist $\varphi_1, \varphi_2 \in D_{A}(x)$ such that $\varphi_1(ax)=\lambda_1$ and $\varphi_2(ax)=\lambda_2$. Let $r \in [0,1]$ and $\varphi=r\varphi_1+(1-r)\varphi_2 \in A^*$. Then $\|\varphi\|= \|r\varphi_1+(1-r)\varphi_2\| \leq |r| \|\varphi_1\|+|(1-r)| \|\varphi_2\| \leq 1$ and $\varphi(x)=(r\varphi_1+(1-r)\varphi_2)(x)=1$. Therefore $\|\varphi\|=1$. So $\varphi \in D_{A}(x)$ and $r\lambda_1+(1-r)\lambda_2 = r\varphi_1(ax)+(1-r)\varphi_2(ax) = \varphi(ax) \in V_{A}(a;x)$.

(3) Since $x$ and $y$ are linearly dependent, there exists $\alpha \in \mathbb C$ such that $y = \alpha x$. Since $\|x\|=\|y\|=1$, we get $|\alpha| =1$. Let $\lambda \in V_{A}(a;x)$. Then there exists $\varphi \in D_A(x)$ such that $\varphi(ax)=\lambda$. Define $\psi(z)=\overline{\alpha}\varphi(z)\, (z \in A)$. Then $\psi \in D_A(y)$, and $\lambda = \varphi(ax) = \varphi(\overline{\alpha}\alpha (ax)) = \overline{\alpha} \varphi(a(\alpha x))= \overline{\alpha} \varphi(ay)= \psi(ay)$. Hence $V_{A}(a;x) \subset V_{A}(a;y)$. Similarly, we can prove that $V_{A}(a;y) \subset V_{A}(a;x)$.

(4) Let $\lambda \in V_B(a;x)$. So there exists $\varphi \in D_B(x)$ such that $\lambda = \varphi(ax)$. By the Hahn-Banach extension theorem, there exists $\widetilde{\varphi} \in A^*$ such that $\widetilde{\varphi}|_B=\varphi$ and $\|\widetilde{\varphi}\|=\|\varphi\|$. Therefore $\widetilde{\varphi} \in D_A(x)$. Hence $\lambda = \varphi(ax) = \widetilde{\varphi}(ax) \in V_A(a;x)$.
Conversely, let $\lambda \in V_A(a;x)$. So there exists $\widetilde{\varphi} \in D_A(x)$ such that $\lambda = \widetilde{\varphi}(ax)$.  Now consider $\varphi =  \widetilde{\varphi}|_B$. Since $x \in S(B)$, we have $\varphi(x) = \widetilde{\varphi}|_B(x) = \widetilde{\varphi}(x)=1$. So $\|\varphi\| = 1$. Thus, $\varphi \in D_B(x)$. Hence $\lambda = \widetilde{\varphi}(ax) = \widetilde{\varphi}|_B(ax) = \varphi(ax) \in V_B(a;x)$.

(5) Let $\lambda \in V_A(a;x)$. There exists $\varphi \in D_A(x)$ such that $\lambda = \varphi(ax)$. Then $\lambda = \varphi(ax) =\lim_{n \rightarrow \infty}\varphi(a_nx) \in \overline{\cup_{n=1}^{\infty}V_A(a_n;x)}$.
\end{proof}

\begin{thm}\label{r2.1}
Let A be any normed algebra. Let $a,b \in A, \, x,y \in S(A)$, and $\alpha \in \mathbb C$. Then
\begin{enumerate}
  \item $V_{A}(\alpha a)=\alpha V_{A}(a)$ and $V_{A}(a+b) \subset V_{A}(a)+ V_{A}(b)$;
  \item $V_A(a)$ is bounded but need not be closed;
  \item If $A$ has a finite dimension, then $V_A(a)$ is compact;
  \item Let $B$ be a subalgebra of $A$. Then $V_B(a) \subset V_A(a) \; (a \in B)$;
  \item Let $a_n \longrightarrow a$ in $A$. Then $\overline{V_A(a)} \subset \overline{\cup_{n=1}^{\infty}V_A(a_n)}$.\label{r2a.1}
\end{enumerate}
\end{thm}

\begin{proof}
(1) This is trivial. Also See \cite[Lemma 2.2]{GaHu:89}.

(2) Note that $|\varphi(ax)| \leq \|ax\| \leq \|a\| \; (x \in S(A), \varphi \in D_A(x))$. Hence $V_A(a)$ is bounded. Let $A = (\ell^1, \|\cdot\|_1)$ with pointwise product. Let $a=\sum_{n=1}^{\infty}\frac{1}{n^2}\delta_n \in A$. Then
 $V_A(a)=(0,1]$. Thus spatial numerical range may not be closed.

(3) The $V_A(a)$ is bounded by Statement (2) above.\\
Let $\lambda \in \overline{V_A(a)}$. Then there exists $\{\lambda_n\} \subset V_A(a)$ such that $\lambda_n$ converges to $\lambda$. So there exist $x_n \in S(A)$ and $\varphi_n \in D_A(x_n)$, where $n \in \mathbb N$ such that $\varphi_n(ax_n) = \lambda_n$. Since $A$ has finite dimension, $S(A)$ is compact. So without loss of generality, we can assume that $x_n \longrightarrow x$ in $S(A)$. Similarly, $S(A^*)$ is compact. So $\{\varphi_n\}$ has a convergent subsequence, say $\varphi_{n_k} \longrightarrow \varphi$. So we can have
\begin{eqnarray*}
  |1-\phi(x)| &=& |\varphi_{n_k}(x_{n_k})-\varphi(x)| \\
   &=& |\varphi_{n_k}(x_{n_k})-\varphi_{n_k}(x)+\varphi_{n_k}(x)-\varphi(x)| \\
   &=& \|\varphi_{n_k}\| \|x_{n_k} - x\| + \|\varphi_{n_k}-\varphi\|\|x\|\\
   & & \longrightarrow 0.
\end{eqnarray*}
Therefore $\varphi(x) = 1$. Thus $\|\varphi\| = 1$ and $\|x\|=1$ and so $\varphi(ax) \in V_A(a)$. Now
\begin{eqnarray*}
 |\lambda_{n_k}- \varphi(ax)| &=& |\varphi_{n_k}(ax_{n_k})-\varphi(ax)|\\
  &=& |\varphi_{n_k}(ax_{n_k})-\varphi_{n_k}(ax)+\varphi_{n_k}(ax)-\varphi(ax)| \\
   &=& \|\varphi_{n_k}\| \|ax_{n_k} - ax\| + \|\varphi_{n_k}-\varphi\|\|ax\|\\
   & & \longrightarrow 0.
\end{eqnarray*}
Thus $\lambda = \lim_{k \rightarrow \infty} \lambda_{n_k}= \varphi(ax) \in V_A(a)$. Hence $V_A(a)$ is closed and so compact.

(4) Let $\lambda \in V_B(a)$. Then $\lambda \in V_B(a;x)$ for some $x \in S(B)$. So there exists $\varphi \in D_B(x)$ such that $\lambda = \varphi(ax)$. By the Hahn-Banach extension theorem, there exists $\widetilde{\varphi} \in A^*$ such that $\widetilde{\varphi}|_B=\varphi$ and $\|\widetilde{\varphi}\|=\|\varphi\|$. Therefore $\widetilde{\varphi} \in D_A(x)$. Hence $\lambda = \varphi(ax) = \widetilde{\varphi}(ax) \in V_A(a;x) \subset V_A(a)$.

(5) Let $\lambda \in V_A(a)$. Choose $x \in S(A)$ and $\varphi \in D_A(x)$ such that $\lambda = \varphi(ax)$. Then $\lambda = \varphi(ax) =\lim_{n \rightarrow \infty}\varphi(a_nx) \in \overline{\cup_{n=1}^{\infty}V_A(a_n)}$.
\end{proof}

The following well-established results follow from Theorems \ref{r2} and \ref{r2.1} above.
\begin{coro}\cite[Section 10]{BoDu:73}\label{cor1}
Let $A$ be a unital normed algebra and let $a \in A$. Then
\begin{enumerate}
\item $V_A(a) = V_A(a; 1)$;
\item $V_A(a)$ is compact and convex;
\item $e^{-1}\|a\| \leq \nu_A(a) \leq \|a\|$;
\item If $B \subseteq A$ is a subalgebra, then $V_B(b) \subseteq V_A(b; 1) \; (b \in B)$;
\item If $B \subseteq A$ is a unital subalgebra with same identity $1$, then $$V_B(b; 1) = V_A(b; 1) \; (b \in B).$$ $\hfill \Box$
\end{enumerate}
\end{coro}

The first statement of the following result is a generalization of \cite[Theorem 2.3]{GaHu:89}; niether the completeness nor the regularity of the norm is required. In the rest of the paper, the set $\overline{co}(K)$ denotes the closed convex hull of $K \subset \mathbb C$.

\begin{thm}\label{Ae}
Let $(A, \|\cdot\|)$ be any non-unital, faithful, normed algebra. Then,
\begin{enumerate}
\item $\overline{co}V_A(a) = V_{A_e^{op}}(a;1)$;
\item $\nu_A(a) = \nu_{A_e^{op}}(a;1)$;
\item $\frac{1}{e}\|a+\lambda1\|_{op} \leq v_{A_e^{op}}(a+\lambda 1; 1) \leq \|a+\lambda1\|_{op} \quad{(a+\lambda1 \in A_e)}$;
\item $\frac{1}{e}\|a\|_{op} \leq v_A(a) \leq \|a\|_{op} \;(a \in A)$;
\item If $\|\cdot\|$ is regular, then $\frac{1}{e}\|a\| \leq v_A(a) \leq \|a\| \;(a \in A)$;
\item If $\|\cdot\|$ is regular and complete, then $0 \in \overline{co}V_A(a) \; (a \in A)$.
\end{enumerate}
\end{thm}

\begin{proof}
(1) For each $a+\lambda 1 \in A_e$, define $L_{a+\lambda1}(x) = ax + \lambda x \; (a \in A)$.
Then $L_{a+\lambda1} \in BL(A)$. We set $\Phi: A_e^{op} \longrightarrow BL(A)$ by $\Phi(a+\lambda1)=L_{a+\lambda 1}$.
Clearly, the map $\Phi$ is an algebra homomorphism. If $\Phi(a+\lambda1)=0$, then for all $x \in A$, $ax + \lambda x=0$ and so $x=\frac{-a}{\lambda} x$ if $\lambda \neq 0$. Thus $\frac{-a}{\lambda}$ for $\lambda \neq 0$, is an identity of $A$, which is ruled out by hypothesis. If $\lambda = 0$, then $ax=0$. By the faithfulness, $a=0$. Therefore $a+\lambda1=0$. Thus $\Phi$ is injective.\\
Also, we have $\|a+\lambda1\|_{op}=\|L_{a+\lambda1}\|$. So $\Phi$ is an isometric algebra isomorphism. Now consider $L_{A_e}=\{L_{a+\lambda1}: a+\lambda1 \in A_e\}$, which is a unital subalgebra of $BL(A)$. Then $\Phi:(A_e, \|\cdot\|_{op}) \longrightarrow (L_{A_e^{op}},\|\cdot\|)$ is an isometric onto algebra isomorphism. This implies that
\begin{eqnarray}\label{eq1}
V_{A_e^{op}}(a+\lambda 1;1) = V_{L_{A_e}}(L_{a+\lambda1};I) \quad{(a+\lambda1 \in A_e)}.
\end{eqnarray}
Hence, by \cite[Theorem 4(ii), page. 84]{BoDu:71}),
\begin{eqnarray*}
\overline{co}V_A(a) &=& V_{BL(A)}(L_a;I)\\
   &=& V_{L_{A_e}}(L_a;I) \quad {(\text{By Corrollary }\ref{cor1}(5))}\\
   &=& V_{A_e^{op}}(a;1)\quad (\text{ By Equation } (\ref{eq1})).
\end{eqnarray*}
Hence $\overline{co}V_A(a) = V_{A_e^{op}}(a;1)$.

(2) This is clear because $\sup\{|\lambda|: \lambda \in \overline{co}(K)\} = \sup\{|\lambda|: \lambda \in K\}$ for any $K \subset \mathbb C$.

(3) This follows from Corollary \ref{cor1}(3).

(4) This is follow from Statement (2) and (3) above.

(5) This is immediate from Statement (4).

(6) Define $\varphi_{\infty}:A_e^{op} \longrightarrow \mathbb C$ as $\varphi_{\infty}(x+\lambda1) = \lambda$. Since $\|\cdot\|$ is complete and regular, $ker\varphi_{\infty} = A$ is closed in $A_e^{op}$. Hence   $\varphi_{\infty}$ is continuous. Therefore we have $0 \in V_{A_e^{op}}(a; 1) \, (a \in A)$. By Statement (1), we get $0 \in \overline{co}V_A(a) \, (a \in A)$.
\end{proof}

The reader should compare Theorem \ref{Ae}(1) with Theorem \ref{Ae1}(1) bellow.

\begin{thm}\label{Ae1}
Let $A$ be any non-unital normed algebra. Then
\begin{enumerate}
\item $\overline{co}(V_A(a) \cup \{0\}) \subset V_{A_e^1}(a;1) \; (a \in A)$;
\item $\frac{1}{e}\|a+\lambda1\|_1 \leq \nu_{A_e^1}(a+\lambda1;1) \leq \|a+\lambda1\|_1 \; (a + \lambda1 \in A_e)$;
\item $\frac{1}{e}\|a\| \leq \nu_{A_e^1}(a;1) \leq \|a\| \; (a \in A)$.
\end{enumerate}

\begin{proof}
(1) Define $\varphi_{\infty}: A_e^1 \longrightarrow \mathbb C$ as $\varphi_{\infty}(x+\lambda1) = \lambda$. Then $\varphi_{\infty} \in D_{A_e^1}(1)$ and so, $0 = \varphi_{\infty}(a) \in V_{A_e^1}(a; 1)$. Hence, by Theorem~\ref{r2.1}(4), $V_A(a) \cup \{0\} \subset V_{A_e^1}(a; 1)$. Since $V_{A_e^1}(a;1)$ is closed and convex as per Theorem 2.1(2), $\overline{co}(V_A(a) \cup \{0\}) \subset V_{A_e^1}(a;1)$.

(2) This is immediate from Corollary \ref{cor1}(3).

(3) Take $\lambda = 0$ in Statement (2).
\end{proof}
\end{thm}

Next result gives all possible relations among $V_A(a), V_{A^{op}}(a), V_{A_e^{op}}(a + \lambda 1)$, and $V_{A_e^1}(a + \lambda 1)$. It could be viewed as a summary of this paper.

\newpage
\begin{thm}
Let $A$ be a faithful normed algebra, let $a \in A$, and let $\lambda \in \mathbb C$. Then
\begin{enumerate}
  \item $\overline{co} V_A(a) = V_{A_e^{op}}(a; 1)$;
  \item $V_{A^{op}}(a) \subseteq \overline{co} V_A(a)$;
  \item $\overline{co} (V_A(a) \cup \{0\}) \subseteq V_{A_e^{1}}(a; 1)$;
  \item $V_{A^{op}}(a) \subseteq V_{A_e^{op}}(a; 1)$;
  \item $V_{A^{op}}(a) \subseteq  V_{A_e^1}(a; 1)$;
  \item $V_{A_e^{op}}(a + \lambda 1; 1) \subseteq V_{A_e^1}(a + \lambda 1; 1)$.
\end{enumerate}
\end{thm}

\begin{proof}
(1) This is Theorem~\ref{Ae}(1) above.

(2) We know that $A^{op}$ is subalgebra of $A_e^{op}$. So by Theorem~\ref{r2.1}(4) and Statement (1) above, we get $V_{A^{op}}(a) \subseteq \overline{co} V_A(a)$.

(3) This is Theorem~\ref{Ae1}(1) above.

(4) $A$ is a subalgebra of $A_e$ and norm on $A$ and $A_e$ are same. Therefore by Theorem \ref{r2.1}(4), we get $V_{A^{op}}(a) \subseteq V_{A_e^{op}}(a; 1)$.

(5) This is follow from Statement (1) and (3) above.

(6) It is true that $V_{A_e^{op}}(a + \lambda1) = V_{A_e^{op}}(a; 1) + \lambda$. So by Statement (1), we get $V_{A_e^{op}}(a + \lambda1; 1) = \overline{co} V_A(a) + \lambda$. Therefore by Statement (3) above, we have $V_{A_e^{op}}(a + \lambda1; 1) \subseteq V_{A_e^1}(a; 1) + \lambda = V_{A_e^1}(a + \lambda1; 1)$.
\end{proof}

\section{Counter Examples}\label{sec:3}

In this section, we intend to give several examples showing that various conditions in our main results cannot be omitted.\\
Let $1 \leq p \leq \infty$. Let $A = (\mathbb C^n, \|\cdot\|_p)$ be a Banach space with pointwise linear operations. Then the dual of $A$ is $A^* =(\mathbb C^n, \|\cdot\|_q)$, where $\frac{1}{p} + \frac{1}{q} = 1$. Therefore, if $\varphi \in A^*$, then there exists unique $y \in A^*$ such that $\varphi = \phi_y$ and $\|\varphi\| = \|y\|_q$, where $\phi_y(x) = \langle x, y \rangle = \sum_{k=1}^{n}x_ky_k \, (x \in A)$.

\begin{exam}
Let $A = (\mathbb C^2, \|\cdot\|_p)$ with the product $xy = (x_1y_1, x_1y_2) = x_1y$. Then $(A, \|\cdot\|_p)$ is a non-unital Banach algebra.  Then, for $a= (a_1, a_2) \in A$ and $x = (x_1, x_2) \in S(A)$,
\begin{eqnarray}
V_A(a; x) = \{a_1\} \text { and } V_A(a) = \{a_1\}.
\end{eqnarray}
\end{exam}

(I) Take $a = (1, 0)$, $a_n = (0,0)$ and $a_n = a \, (n \geq 2)$. Then clearly $a_n \longrightarrow a$. So we have $\overline{V_A(a;x)} = \{1\}$ and $\overline{\cup_{n=1}^{\infty}V_A(a_n;x)} = \{1, 0\}$. Hence the inclusion in Theorems \ref{r2}(5) and \ref{r2.1}(5) are proper for the sequence $\{a_n\}$.

(II) Let $A = (\mathbb C^2, \|\cdot\|_{\infty})$. Let $a=(1,1), \, x=(1,0) \text{ and } y=(0,1)$. Then $x$ and $y$ are linearly independent, but $V_A(a,x) = V_A(a,y) = \{1\}$. So, the converse of Theorem~\ref{r2}(3) is not true.

(III) Let $A = (\mathbb C^2, \|\cdot\|_1)$. Then $A_e^1 \cong \mathbb C^3$ is a Banach algebra with the product $xy = (x_1y_1+x_1y_3+x_3y_1, x_1y_2 + x_2y_3 + x_3y_2, x_3y_3)$. The dual of $A_e^1$ is $(\mathbb C^3, \|\cdot\|_{\infty})$. Let $\varphi \in (A_e^1)^*$ and let $a \in A_e^1$. Then there exists $y=(y_1, y_2, y_3) \in \mathbb C^3$ such that $\varphi = \phi_y$ and $\|\varphi\| = \|y\|$. Now
\begin{eqnarray*}
   V_{A_e^1}(a;1) & = & \{\varphi(a):\varphi \in D_{A_e^1}(1)\}\\
    & = &  \{\phi_y(a):y \in \mathbb C^3 \text{ and } \max\{|y_1|,|y_2|,|y_3|\} = 1 = y_3\}\\
    & = &  \{a_1y_1 + a_2y_2:\max\{|y_1|, |y_2|\} \leq 1\}\\
    & = &  \{\|a\|_1z: |z| \leq 1\}.
 \end{eqnarray*}
In particular, taking $a = (1,0)$, we get $\overline{co}(V_A(a) \cup \{0\}) = [0,1]$ and $V_{A_e^1}(a) = \{z \in \mathbb C: |z| \leq 1\}$. Thus the inclusion in Theorem \ref{Ae1}(1) may happen to be proper.

(IV) Let $A = (\mathbb C^2, \|\cdot\|_1)$. Then $\|\cdot\|_1$ is not regular. Let $a = (1,0) \in A$. Since $V_A(a)=\{a_1\}$ and so $0 \notin \overline{V_A(a)}$. Hence Theorem \ref{Ae}(6) is not true if the norm is not regular.

\begin{exam}
Let $A = (\mathbb C^2, \|\cdot\|_p)$ with the product $xy = (x_1y_1, x_2y_1) = xy_1$. Then $(A, \|\cdot\|_p)$ is a non-unital Banach algebra and $\|\cdot\|_p$ is a regular norm on $A$.
\end{exam}

(I) Let $A = (\mathbb C^2, \|\cdot\|_{\infty})$. Let $x \in S(A)$ and $\varphi = \phi_y \in D_A(x)$. Then we get $ |y_1|+|y_2| = 1 = x_1y_1 + x_2y_2$ and $\varphi(ax) = \phi_y(ax_1) = a_1x_1y_1 + a_2x_1y_2 = a_1r + a_2(1-r)e^{i\theta}$ for some $r \in [0,1]$ and $\theta \in (-\pi, \pi]$.
Since $x \in S(A)$ is arbitrary,
$$V_A(a) \subset \{a_1r + a_2(1-r)e^{i\theta} :r \in [0,1]  \text{ and } \theta \in (-\pi, \pi]\}.$$
Conversely, let $r \in [0,1]$ and $\theta \in (-\pi, \pi]$. Take $x = (e^{-i\theta}, 1)$ and $y = (r e^{i\theta}, 1-r)$. Then $\|x\|_{\infty}=1, \, \phi_y \in D_A(x)$, and $a_1r+a_2(1-r)e^{i\theta} = \phi_y(ax) \in V_A(a;x) \subset V_A(a)$.
Thus $V_A(a) = \{a_1r + a_2(1-r)e^{i\theta} : r \in [0,1], \text{ and } \theta \in (-\pi, \pi]\}$.

(II) Let $A = (\mathbb C^2, \|\cdot\|_1)$. Then, as per the similar arguments in (I) above,
$$V_A(a) = \{a_1r + a_2re^{i\theta} : r \in [0,1], \text{ and } \theta \in (-\pi, \pi]\}.$$

(III) It follows from (I) and (II) that $V_A(a)$ is different in $(\mathbb C^2, \|\cdot\|_{\infty})$ and $(\mathbb C^2, \|\cdot\|_1)$ even though the two norms are equivalent.

(IV) Consider the subalgebra $B = \mathbb C \times \{0\}$ of $A$ and take $a = (1,0) \in B$ which is the identity of $B$ and so $V_B(a) = \{1\}$. But $V_A(a) = [0,1]$ as (II) above. So $V_B(a) \subsetneq V_A(a)$. Thus the set inclusion in Theorem 2.2(4) may happen to be proper.

\begin{exam}
Let $A$ be a unital Banach algebra with identity $1$ such that $A \neq \{0\}$. Take $a = 1$. Then $V_A(1;1) = \{1\}$. So, $0 \notin \overline{co}V_A(1; 1)$. Thus Theorem \ref{Ae}(6) is true only for non-unital case.
\end{exam}

\begin{exam}
Let $A = (\ell^1, \|\cdot\|_1)$ with pointwise operations. Then $\|\cdot\|_1$ is not regular. Define $f_n(k) = 1 \, (k \leq n)$ and $f_n(k) = 0 \, (k > n)$. Then, by Theorem \ref{Ae}(4), $\nu_A(f_n) \leq \|f_n\|_{op} = \|f_n\|_{\infty} = 1$. But $\|f_n\|_1 = n \; (n \in \mathbb N)$. Hence Theorem \ref{Ae}(5) is not true, if the norm $A$ is not regular.
\end{exam}

\noindent
\textbf{Acknowledgement}: The second author is thankful to Council of Scientific and Industrial Research (CSIR), New Delhi, India, for providing Senior Research Fellowship.


\begin{thebibliography}{9}
\bibitem{ArMu:07} J. Arhippainen, and V. Muller, \emph{Norms on unitizations of Banach algebras revisited}, Acta Mathematica Hungarica, 114(3)(2007)201-204.
\vspace{0.2cm}
\bibitem{BoDu:71} F. F. Bonsall, and J. Duncan, \emph{Numerical ranges of operators on normed spaces and of elements of normed algebras}, Cambrige University Press, 1971.
\vspace{0.2cm}
\bibitem{BoDu:73} F. F. Bonsall, and J. Duncan, \emph{Complete normed algebras}, Springer, Berlin, 1973.

\vspace{0.2cm}
\bibitem{DePa:22} H. V. Dedania, and A. B. Patel, \emph{The spectral extension property in the unitization of Banach algebras}, Bull. of Calcutta Math. Soc., 114(5)(2022)759-766.

\vspace{0.2cm}
\bibitem{DePa:23} H. V. Dedania, and A. B. Patel, \emph{On the convexity of the spatial numerical range in non-unital normed algebras}, Communicated.
\vspace{0.2cm}
\bibitem{GaHu:89} A. K. Gaur, and T. Husain, \emph{Spatial numerical ranges of elements of Banach algebras}, Internat. J. Math. \& Math. Sci., 12(4)(1989)633-640.
\vspace{0.2cm}
\bibitem{GaKo:91} A. K. Gaur, and Z. V. Kovarik, \emph{Norms, states and numerical ranges on direct sums}, Analysis, 11(1991)155-164.
\vspace{0.2cm}
\bibitem{GaKo:93} A. K. Gaur, and Z. V. Kovarik, \emph{Norms on unitizations of Banach algebras}, Proc. Amer. Math. Soc., 117(1)(1993)111-113.
\end{thebibliography}
\end{document}